\documentclass{amsart}

\usepackage[latin1]{inputenc}
\usepackage{amsfonts,amsmath}
\usepackage{amsthm}
\usepackage{enumerate}
\date{\today}
\newtheorem{teo}{Theorem}

\newtheorem{lemma}{Lemma}

\theoremstyle{remark}

\newcommand{\RR}{\mathbb{R}}
\newcommand{\ZZ}{\mathbb{Z}}

\newcommand{\NN}{\mathbb{N}_0}

\newcommand{\cqd}{\hfill$\Box$}

\title[The Hermite-Schr\"odinger 
equation]{On the boundary convergence  of solutions to the 
Hermite-Schr\"odinger equation}

\subjclass[2000]{Primary 42B35, Secondary 42C10, 35K15}
\keywords{Hermite expansions, Schr\"odinger equation, Strichartz estimates, 
boundary convergence}

\author{Peter Sj\"ogren}

\author{J.L. Torrea}

\address[P. Sj\"ogren]{Mathematical Sciences, University of Gothenburg, 
SE-412 96 G\"oteborg, Sweden and Mathematical Sciences,
Chalmers, SE-412 96 G\"oteborg, Sweden}
\email{peters@chalmers.se}

\address[J. L. Torrea]{Departamento de Matem\'aticas and ICMAT  CSIC-UAM-UCM-UC3M, Facultad de Ciencias, Universidad Aut\'onoma de Madrid, 28049 Madrid, Spain}
\email{joseluis.torrea@uam.es} 

\thanks{The authors were partially supported by the Ministerio de Educaci\'on y Ciencia (Spain),
through grant SAB-2006-0092 and grant MTM2008-06621-C02-01, respectively.}

\begin{document}


\begin{abstract}
  In the half-space $\RR^d \times \RR_+$, we consider the  
Hermite-Schr\"odinger equation  $i\partial u/\partial t = -
\Delta u + |x|^2 u$, with given boundary values on  $\RR^d$.
 We prove
a formula that links the solution of this
problem to that of the classical Schr\"odinger equation.
It shows that mixed norm estimates for  the 
Hermite-Schr\"odinger equation 
can  be obtained immediately from
those known in the classical case. In one space dimension, we  deduce
sharp  pointwise convergence results at the boundary, by means of this
link.
\end{abstract}

\maketitle

\dedicatory{The authors dedicate this paper to the memory of  Andrzej Hulanicki.
Both of us knew Andrzej since the 1970's. Since then he has been like an invariant for us. We have enjoyed the high quality of his mathematics, his capacity of work, his ability to organize important mathematical events, his generosity when sharing ideas and his sympathy. All this, and even his age, seemed to be invariant during these decades.}

\section{Introduction and results}

  The solution of the classical, free Schr\"odinger  equation
in the half-space $\RR^d \times \RR_+$ with variables $(x, t)$,
$$
\left\{
\begin{array}{ll}
  i\frac{\partial u}{\partial t}
&
 =
-\Delta u
\\[0.3cm]
u(\cdot,0)
&
 =
f
\end{array}
\right.
$$
can  be written $u(x,t) = e^{it\Delta}f(x)$, for $f \in L^2(\RR^d)$.
For $p, q \in [1, \infty]$, one measures the size of $u$ by means of
the mixed norm 
\[
\|u\|_{L^q_t(\mathbb{R},L^p_x(\mathbb{R}^d))} = 
\left(\int_{\mathbb{R}} \left( \int_{\mathbb{R}^d}|u(x, t)|^p\,dx\right)^{q/p}\,dt\right)^{1/q},
\]
with the obvious interpretation for  $p$ or $q=\infty$. The Strichartz 
estimate
\begin{eqnarray}\label{stD}
\|u\|_{L^q_t(\mathbb{R},L^p_x(\mathbb{R}^d))} \le C_{d,p} \|f\|_{L^2(\mathbb{R}^d)}, \qquad 
\end{eqnarray}
is known to hold if and only if 
\begin{equation}\label{assump}
 \frac{d}{p}+\frac2{q}= \frac{d}{2} \qquad \mathrm{and}\qquad
\left\{
\begin{array}{lll}
  2 \le p \le \infty &\mathrm{for} \;\; d=1\\
2 \le p < \infty &\mathrm{for}\;\; d=2\\
 2 \le p \le 2d/(d-2) &\mathrm{for} \;\;  d \ge 3.
\end{array}
\right.
\end{equation}
This is due essentially to J.~Ginibre and G.~Velo \cite{GV}.  M.~Keel 
and T.~Tao  \cite{KT} obtained the endpoint results.

Results about the pointwise convergence of $u(x,t)$ at the boundary  
are also known, for $f$ in  Sobolev spaces
\[
 W^s(\mathbb{R}^d) = \{f \in L^2(\mathbb{R}^d):
(I-\Delta)^{s/2}f \in L^2\}.
\]
For $d=1$,   L. Carleson \cite{Ca} and 
B. Dahlberg and C. Kenig        \cite{DK} have proved that
$e^{it\Delta}f \to f$ a.e.  as 
$t \to 0^+$  for all $f \in  W^s(\mathbb{R})$ if and only if
$s \ge 1/4$.


In this paper, we consider the same questions for the Hermite operator
 $$H = -\Delta + |x|^2, \qquad x \in \mathbb{R}^d. $$
Thus $u$ will be the solution $u(x,t) = e^{-itH}f(x)$ to the Hermite-Schr\"odinger 
equation in  $\RR^d \times \RR_+$   with given boundary values,
\begin{eqnarray}
  \label{pr}
  \left\{
\begin{array}{ll}
 i\frac{\partial u}{\partial t}
&
 =
H u
\\[0.3cm]
u(\cdot,0)
&
 =
f.
\end{array}
\right.
\end{eqnarray}   
As in the classical case, the  Strichartz estimate
\begin{eqnarray}\label{stH}
\|e^{-itH}f\|_{L^q_t((0,2\pi),L^p_x(\mathbb{R}^d))} \le C_{d,p,q} \|f\|_{L^2(\mathbb{R}^d)}
\end{eqnarray}
holds under the  assumption (\ref{assump}); see H.~Koch and 
D.~Tataru \cite{KoTa}. 
 Moreover, since the interval of integration in the $t$ variable is now
bounded, (\ref{stH}) remains true if the equality in (\ref{assump}) is 
replaced by
the inequality $d/p + 2/q \ge d/2$.

Our Lemma \ref{fund}
in Section 2 gives an explicit relation
between the two solution operators $e^{-itH}$ and $e^{it\Delta}$. It makes
it easy to prove the following result,  which  implies that the
estimates (\ref{stD}) and (\ref{stH}) are actually equivalent when the
equality in (\ref{assump}) holds.

\begin{teo}\label{T1} Let $1 \le p, q \le \infty,$ and 
assume that $\displaystyle \frac{d}{p}+\frac2{q} =\frac{d}{2}.$   Then
for $f \in L^2$
$$
 \|e^{-itH}f\|_{L^q_t((0,\pi/4),L^p_x(\mathbb{R}^d))} =
 \|e^{it\Delta}f\|_{L^q_t((0,\infty),L^p_x(\mathbb{R}^d))}. 
$$
\end{teo}

As we shall see below, it does not matter whether the $t$ interval 
in (\ref{stH})
is $(0,2\pi)$ or $(0,\pi/4)$;  the two mixed norms obtained are proportional
for real functions $f$. 

In the case $d=1$, we shall also consider the almost everywhere convergence 
as $t\rightarrow 0^+$  of the solution $e^{-itH}f$, to the initial data.  To 
state these results, we  use both $ W^s(\mathbb{R})$ and the Sobolev 
spaces associated to 
 $H,$ defined by 
\[
 W_H^s(\mathbb{R}) = \{f \in L^2(\mathbb{R}):
H^{s/2}f \in L^2\} 
\]
with the obvious norm. These spaces have
been introduced by S. Thangavelu \cite{Th}. We point out that there is
a continuous inclusion
$W^s_H \subset W^s$, see B. Bongioanni and J.L. Torrea 
\cite[Theorem 3(i)]{BT}. 

Yajima \cite{Ya} proved the a.e. convergence   $e^{-itH}f \to f$  
as $t\rightarrow 0^+$ for $f$ in the intersection 
$W^s(\mathbb{R}) \cap L^1$, with $s>1/2$. Then
 Bongioanni and  Rogers
\cite{BR} obtained the same convergence  for 
$f \in W_H^s(\mathbb{R})$, with $s > 1/3$. The following result is sharp
with respect to both types of Sobolev spaces.

 \begin{teo} \label{convergence} Let $d=1$.\\
(i) Assume 
$f \in W^{1/4}(\mathbb{R})$. Then for a.a. $x \in \mathbb{R}$  the 
function  $t \mapsto e^{-itH}f(x), \; 0<t<\pi/8$, will, after modification
on a null set, be continuous with limit $f(x)$ as $t \to 0^+$.\\
(ii) If  $s < 1/4$, there exists an 
$f \in W^s(\mathbb{R})$ such that for all $x$ in a set of positive
measure the function  $t \mapsto e^{-itH}f(x)$ does not converge as
 $t \to 0^+$, even after any modification on a null set. \\
(iii) The statements in  (i) and (ii) also hold if the spaces 
$W^s(\mathbb{R})$ are  replaced by 
$W_H^s(\mathbb{R})$.
\end{teo}

Bongioanni and  Rogers obtained their convergence result via a sharp
global maximal operator estimate from 
$W_H^s$ into $L^p$. The relevant maximal function is
\[
\mathcal{M}f(x) = \underset{0<t<\pi/8}{\mbox{ess sup}} \; |e^{-itH}f(x)|.
\]
The proof of Theorem \ref{convergence}(i) is based on a local $L^1$
estimate for $\mathcal{M}$.
The following result  says that there is no 
global $L^p$ estimate for  $\mathcal{M}$ from $W^s$ into $L^p$.

\begin{teo} \label{noway} Let $d=1$ and $s>0$. The operator $\mathcal{M}$
does not map $W^s(\mathbb{R})$ boundedly into $L^p$ nor 
into weak $L^p$ for any $p< \infty$, 
and $\mathcal{M}$ maps  $W^s(\mathbb{R})$   boundedly into  
$L^\infty$ if and only if  $s > 1/2.$
\end{teo}

By $c>0$ and $C<\infty$ we denote many different constants.

\section{ Some key formulas;  proof of  Theorem \ref{T1} }
Let $h_n(x),\: n \in \NN,$ denote the Hermite functions in $\RR$, normalized
in $L^2$. By $\Phi_\mu, \: \mu \in \NN^d$, we denote the $d$-dimensional,
normalized
Hermite functions, which are simply the tensor products of the $h_n$.
 See further Thangavelu \cite[Sect. 1.1]{ThangaLect}.

The semigroup $e^{-tH},\:t>0$, generated by $H$ can be defined also with a
complex parameter $z$ instead of $t$, for  $\Re z > 0$. Moreover, for these
$z$ the operator  $e^{-zH}$ is given by integration against the kernel
\begin{equation}\label{sum}
K_z(x,y)  = \sum_{\mu \in \NN^d} e^{-(2|\mu|+d)z} \Phi_\mu(x)\Phi_\mu(y). 
\end{equation}
For real and for complex parameter values, this series can be summed.
The sum is the well-known Mehler kernel, which can be found  for instance  in 
\cite[equation (4.1.3) p. 85]{ThangaLect}. For  $\Re z > 0$ one has
\begin{multline*}
  K_{z}(x,y) =  \frac1{(2\pi \sinh 2z)^{d/2}}
  \exp{\left(\frac{1}2\Big( -\coth 2z \: (|x|^2 + |y|^2) + \frac2{\sinh 2z}  \: x\cdot y  \Big)\right)}.
\end{multline*}
This expression is well defined also for $z$ on the imaginary axis,
except at the multiples of $i\pi/2$. Indeed, for 
$t \in \RR \setminus \frac{\pi}2 \ZZ$ we get 
\begin{eqnarray} \label{definitva}
  K_{it}(x,y) =  \frac{e^{-i\pi d/4} }{(2\pi \sin 2t)^{d/2}}   \notag
  \exp{\left(\frac{i}2\Big( \cot 2t \, (|x|^2 + |y|^2) - \frac2{\sin 2t}  \, x\cdot y  \Big)\right)}\\     
=  \frac{e^{-i\pi d/4} }{(2\pi \sin 2t)^{d/2}}
\exp{\left(\frac{i}2\Big(\cot 2t \, 
\Big|y-\frac{x}{\cos 2t}\Big|^2 -  \tan 2t \, |x|^2\Big)\right)}.
\end{eqnarray}
By analytic continuation from $\Re z > 0$, one sees that the 
argument of the quantity  $(2\pi \sin 2t)^{d/2}$     occurring here
should be chosen as $[2t/\pi]\pi d/2$.
 One can also check that
integration against this
kernel gives the solution of the problem (\ref{pr}), at least for test
 functions $f$. Since $K_{it}$ is the kernel of $e^{-itH}$, we shall often write
$K_{it}f$ instead of $e^{-itH}f$. Clearly, each operator   $e^{-itH}$    is bounded on   $L^2$. 

The Hermite functions $h_n$ are  real-valued and have the same parity
as the index $n$. From (\ref{sum}), it 
follows that $K_{\overline{z}}(x,y) = \overline{K_z(x,y)}$, and also that
$K_{z+i\pi/2}(x,y) = e^{-i\pi d/2} K_z(-x,y)$. Here $\Re{z} > 0$, but if 
$t \in \RR$ is not a multiple of $\pi/2$, we also conclude that
\[
K_{-it}(x,y) = \overline{K_{it}(x,y)} \qquad \mathrm{and} \qquad
K_{i(t+\pi/2)}(x,y) = e^{-i\pi d/2} K_{it}(-x,y).
\]
For real functions $f$, it follows that the $L^p(\RR^n)$ norm of
 $e^{-itH}f$ is even and $\pi/2$-periodic as a function of  $t$, and thus determined 
by its values for $0<t<\pi/4$.

We shall compare the operators   $e^{-itH}$   and  $e^{it\Delta}$  by finding a 
link between their kernels. The kernel of  $e^{it\Delta}$    is the 
standard Schr\"odinger kernel

\[
L_{it}(x,y) = e^{-i\pi d/4} \frac1{(4\pi t)^{d/2}} 
\exp{\Big(i\frac1{4t} |y-x|^2\Big)}.
\]
Instead of $e^{it\Delta}f$, we shall often write $L_{it}f$.

\begin{lemma}  \label{fund}
For any $f \in L^2$ and any $v>0$,
 \begin{equation*}
K_{i\frac{\arctan v}2}f(x) =  \exp{(-iv |x|^2/2)} \, (1+v^2)^{d/4}\,
L_{iv/2}f(x\,\sqrt{1+v^2}\,).
  \end{equation*} 
\end{lemma}

\begin{proof}
  For $0 < t < \pi/4$,  we   let $\tan 2t = v$ in (\ref{definitva}) and get 
\begin{eqnarray*}
K_{i\frac{\arctan v}2}(x,y)&\\
=& e^{-i\pi d/4} \left(\frac{\sqrt{1+v^2}}{2\pi v}\right)^{d/2} \exp{\Big(-i\frac{v}{2} |x|^2\Big)}\,
\exp{\Big(i\frac1{2v} \, \Big|y-x\,\sqrt{1+v^2}\Big|^2\Big)} \\
=& \exp{(-iv |x|^2/2)}  (1+v^2)^{d/4}
L_{iv/2}(x\,\sqrt{1+v^2}, y).
\end{eqnarray*}
Integrating against $f(y)\,dy$, we obtain the desired equation when
 $f \in C_0^\infty$. The general case then follows  by continuity in $L^2$.
\end{proof}

\noindent \textit{Proof of Theorem \ref{T1}.} \hskip1mm Assuming 
$p,q < \infty,$ we get
\begin{multline*}
\int_0^{\pi/4} \Big(\int_{\RR^d} |K_{it} f(x)|^p\,dx\Big)^{q/p}\,dt =
\int_0^{\infty} \Big(\int_{\RR^d} |K_{i\frac{\arctan v}2} f(x)|^p\,dx\Big)^{q/p}
\frac12\frac1{1+v^2}\,dv \\
=\frac12 \int_0^{\infty} \Big(\int_{\RR^d} |(1+v^2)^{d/4}
L_{iv/2}f(x\,\sqrt{1+v^2})|^p\,dx\Big)^{q/p}
\frac1{1+v^2}\,dv \\
=\frac12 \int_0^{\infty} \Big(\int_{\RR^d} |
L_{iv/2}f(x)|^p\,dx\Big)^{q/p}
(1+v^2)^{-q(d/p+2/q-d/2)/2}\,dv \\
=  \int_0^{\infty} \Big(\int_{\RR^d} |
L_{iv}f(x)|^p\,dx\Big)^{q/p}\,dv. 
\end{multline*}
The cases when  $p$ or $q$ is infinite are similar.
\cqd

\section{ Proof of Theorem \ref{convergence}    
}

 From now on, $d=1$.
In this section, we shall need the following  estimate, which is based on
Carleson's lemma in \cite[p. 24]{Ca}. It can also be seen as a limit case
of a lemma due to Kenig and A. Ruiz \cite[Lemma 2]{KR} (cf. (7) below), 
but we prefer to give a direct proof.

\begin{lemma}\label{Peter}  Let $a$ and $b$ be real numbers with 
$(a,b) \neq (0,0).$  Then for any interval $J$
$$\Big| \int_J e^{i(a t +b t^2)}\, \frac{d t}{|t|^{\frac12}} \Big| \le C \min( |a|^{-\frac12},|b|^{-\frac14}),$$
where  $C$ is an absolute constant. If $J$ is unbounded,  the integral here is 
interpreted as the limit of the integrals over bounded intervals increasing
to $J$.
\end{lemma}

\begin{proof} Assume first $b=0$. By homogeneity, we need then only 
consider the case $a=1$, which is easy.

When $b \ne 0$, we see 
by taking the
conjugate that we may  assume $b>0.$ Let $u= b^{\frac12} t$ and $A = -a b^{-\frac12}/{2}.
$ Then 
\begin{multline*} \int_J e^{i(a t +b t^2)}\, \frac{d t}{|t|^{\frac12}}  
= b^{-\frac14} \int_{J'} e^{i(-2Au+u^2)} \frac{du}{|u|^{\frac12}} 
= e^{-iA^2} b^{-\frac14}\int_{J'} e^{i(u-A)^2} \frac{du}{|u|^{\frac12}},
\end{multline*}
for some interval $J'$.
The lemma is equivalent to   the following  claim:
\begin{equation}\label{claim} \Big| \int_{J'} e^{i(t-A)^2} \frac{dt}{|t|^{\frac12}}\Big|\,\, \le \,\,C\min(1, |A|^{-\frac12 }). \end{equation}

Without loss of generality, we may assume $A\ge 0$. Consider first the case
$0\le A \le 2$. Then we split the integral in (\ref{claim}) and integrate by parts in the second term, getting
\begin{multline*}
 \Big|\int_{J'} e^{i(t-A)^2} \frac{dt}{|t|^{\frac12}}\Big| \le \Big| \int_{|t| <4}
\chi_{J'} e^{i(t-A)^2} \frac{dt}{|t|^{\frac12}} \Big|+\Big|\int_{|t|>4}\chi_{J'}e^{i(t-A)^2} \frac{dt}{|t|^{\frac12}}\Big|  \\ \le  \int_{|t| <4} \frac{dt}{|t|^{\frac12}} +|\mathrm{integrated\; terms}| + \Big| \int_{|t| >4}\chi_{J'} e^{i(t-A)^2} \frac{d}{dt}\Big(\frac1{2(t-A)|t|^{\frac12}}\Big) dt\Big| \\
 \le  C + C + C\int_{|t|>4} \frac{dt}{|t|^{\frac52}} \le C.
 \end{multline*}

Now let $A>2$. We begin by observing that
\[
 \Big| \int_{|t| < 1/{A}} \chi_{J'} e^{i(t-A)^2} \frac{dt}{|t|^{\frac12}}\Big| \le  \int_{|t| < 1/{A}} \frac{dt}{|t|^{\frac12}} \le C A^{-\frac12}
\]
and
\[
\int_{|t-A|< 1  }\chi_{J'} \frac{dt}{|t|^{\frac12}} \le C A^{-\frac12}.
\] 
In the remaining integral,  taken over the set 
${\{t \in J': |t| > 1/{A} \hbox{ and } |t-A|> 1 \}}$, we integrate by parts, as
above. The integrated terms will then be controlled
by the values of $(t-A)^{-1}|t|^{-1/2}$ at a few points in the set
$\{|t| \ge 1/A \hbox{ and } |t-A| \ge 1\}$,
and those values are all bounded by $CA^{-1/2}$. 
Thus we need only consider the integral
\begin{multline*}
 \Big|\int_{\{t \in J': |t| > 1/{A} \hbox{ and } |t-A|> 1\} } e^{i(t-A)^2} \frac{d}{dt}\Big(\frac{1}{ (t-A)|t|^{\frac12}}\Big) dt \Big| \\ \le   \int_{\{|t| > 1/{A} \hbox{ and } |t-A|> 1\} } \frac{1}{ (t-A)^2|t|^{\frac12}} dt  +   \int_{\{|t| >1/{A} \hbox{ and } |t-A|> 1 \}} \frac{1}{ |t-A||t|^{\frac32}} dt \\ = I +II.  \quad \quad \quad \quad \quad 
\end{multline*}
We split each of the integrals $I$ and $II$ thus defined into parts given
by $|t|<A/2$ and $|t|>A/2$. For $I$ we get
\begin{multline*} I = \int_{\{1/A < |t| < A/2\}}+ \int_{\{|t-A|>1 \,\, \hbox{and}\,\, |t| > A/2\}} \frac{1}{ (t-A)^2|t|^{\frac12}} dt\\ \le 
C\int_{|t|<A/2} \frac{1}{ A^2|t|^{\frac12}}dt+ C\int_{|t-A|>1} \frac{1}{ (t-A)^2A^{\frac12}}\,dt
\le CA^{-\frac12}.
\end{multline*}
 and similarly
\begin{multline*} II = \int_{\{1/A < |t| < A/2\}}+ \int_{\{|t-A|> 1 \,\, \hbox{and}\,\, |t| > A/2\}} \frac{1}{ |t-A||t|^{\frac32}} dt\\  \le 
C\int_{|t|>1/{A}} \frac{1}{A |t|^{\frac32}}dt+ \int_{|t|>{A}/{2}} \frac{1}{ |t|^{\frac32}}dt
\le CA^{-\frac12}.
\end{multline*}

The claim is verified, and Lemma \ref{Peter} 
is proved.
\end{proof}

The  maximal function estimate in the next lemma will enable us to prove 
 Theorem~\ref{convergence}(i). For $f \in C_0^\infty$, the function
 $e^{-itH}f(x) = K_{it}f(x)$ is continuous in
  $(x,t) \in \RR \times \overline{\RR_+}$ if defined as $f(x)$ for
$t=0$, as easily
verified with Fourier transforms. In the definition of $\mathcal{M}f$, 
one can for  $f\in C_0^\infty$
obviously replace the essential supremum by an ordinary supremum.

\begin{lemma}\label{essmaxest}
Let $I$ be a bounded interval.   
Then for any   $f \in C_0^\infty(\RR)$,
\begin{equation} \label{maxK}
   \int_I   
\mathcal{M}f(x)\,dx \le C \|f\|_{W^{1/4}},  \qquad C = C(I).
\end{equation}
 \end{lemma}

Before proving this lemma, we use it to prove Theorem~\ref{convergence}(i).
Given $f \in W^{1/4}$,  we take a sequence  $f_j  \in C_0^\infty, \;
j=1, 2, ...,$ with 
$\|f_j - f\|_{W^{1/4}} < 2^{-j}$. Applying Lemma \ref{essmaxest}
 to $f_j - f_{j+1}$,
whose $W^{1/4}$ norm is less than $2^{1-j}$, we get 
\begin{equation}
  \label{eq:split}
\int_I \sup_{0<t<\pi/8}|K_{it}f_j(x) - K_{it}f_{j+1}(x)|\,dx \le C2^{-j}. 
\end{equation}
Here the supremum can be taken over $0 \le t <\pi/8$, 
since each function $K_{it}f_j(x)$ is continuous in $\RR \times 
[0, \pi/8)$
with  the value $f_j(x)$ at $(x,0)$.
The integrals in  (\ref{eq:split}) have a finite sum in $j$, so that
\begin{equation*}
  \sum_{j=1}^\infty \, \sup_{0 \le t <\pi/8}|K_{it}f_j(x) - K_{it}f_{j+1}(x)|
\end{equation*}
is finite for a.a.\ $x\in I$. But for any fixed $x$ with this property,
the functions $t\mapsto K_{it}f_j(x)$ will converge, uniformly in 
$0 \le t <\pi/8$, to a continuous function $u_x(t).$ On the other hand,
 $K_{it}f_j(x) \to K_{it}f(x)$ in $L^2(I\times (0,\pi/8))$,
and $f_j \to f$ in $L^2(\RR)$. We conclude that for a.a.\ $x$, the
function $t\mapsto K_{it}f_j(x)$ must coincide with the continuous
function  $u_x(t)$ for a.a.\ $t \in (0,\pi/8)$ and, moreover, $u_x(0) = f(x)$.
This implies Theorem~\ref{convergence}(i). 

\vskip2mm

\noindent \textit{Proof of Lemma  \ref{essmaxest}.} \hskip2mm
Because of Lemma \ref{fund},  one can replace $\mathcal{M}f(x)$ by
\[
\sup_{0<v<1}\; \left|L_{iv/2}f(x\,\sqrt{1+v^2}\,)\right|
\]
when proving (\ref{maxK}). 
 It is clearly enough to
show that  for all  $f\in C_0^\infty$ 
\begin{equation}
  \label{essmax}
  \int_I \,
\sup_{0<v<1}\;
\Re_+ L_{iv/2}f(x\,\sqrt{1+v^2}\,)\,dx \le C \|f\|_{W^{1/4}},
\end{equation}
where $\Re_+$ denotes the positive part of the real part.

 We first compare  the 
integrals over $I$ of 
\[\sup \,\Re_+ L_{iv/2}f(x\,\sqrt{1+v^2}\,) \qquad \mathrm{and} \qquad
 \sup \,\Re\, L_{iv/2}f(x\,\sqrt{1+v^2}\,),
\]
 where both suprema are taken 
over $0<v<1$.
They differ only on the set 
$M = \{x \in I: \sup\, \Re\, L_{iv/2}f(x\,\sqrt{1+v^2}\,)<0\}.$ 
Since $L_{iv/2}f(x\,\sqrt{1+v^2}\,)$
converges pointwise to $f$ as $v \to 0^+$, we have 
$\sup \Re L_{iv/2}f(x\,\sqrt{1+v^2}\,) \ge \Re f(x)$ for all $x$, and so
\begin{eqnarray*}
&&\int_I \sup\, \Re_+ L_{iv/2}f(x\,\sqrt{1+v^2}\,)\, dx \\  &=&
\int_I \sup\, \Re\, L_{iv/2}f(x\,\sqrt{1+v^2}\,) \,dx 
- \int_M \sup\, \Re\, L_{iv/2}f(x\,\sqrt{1+v^2}\,) \, dx 
\\ &\le &
 \int_I \sup\, \Re\, L_{iv/2}f(x\,\sqrt{1+v^2}\,) \,dx + \int_M (-\Re f(x))\,dx \\
&\le & \int_I \sup\, \Re \,L_{iv/2}f(x\,\sqrt{1+v^2}\,) \,dx +  C \|f\|_{W^{1/4}};
\end{eqnarray*}
here the last step went via an $L^2$
estimate.

This means that we can replace  $\Re_+$ by $\Re$ when we prove 
(\ref{essmax}) for   $f \in C_0^\infty$. We shall use the
 method  of Kolmogorov-Seliverstov-Plessner, see also
 Carleson \cite[Theorem, p.~24]{Ca}.
It is enough to let  $v=v(x)$ be a measurable  function of 
$x \in I$ with $0<v(x)<1$ and to prove that
$$\Re \int_I \ L_{iv(x)/2}f(x\,\sqrt{1+v(x)^2}\,)\, dx   \le C \|f\|_{W^{1/4}},$$
with $C = C(I)$ independent of $v(x)$ and $f.$ 

 We define  the Fourier transform  by 
$\displaystyle \hat{h}(\xi) =
 \int_{\mathbb{R}} h(x) e^{-i x \xi }\, dx$ and observe that 
$\widehat{L_{it}}(\xi) = \exp(-it|\xi|^2)$. This leads to
\begin{eqnarray*}
2\pi \lefteqn{\Big|\int_I  L_{iv(x)/2}f(\,x\,\sqrt{1+v(x)^2}\,)  \, dx \Big| =\Big|\int_{-\infty}^{\infty} \hat{f}(\xi)\int_I e^{ix\,\xi\,\sqrt{1+v(x)^2}} e^{-iv(x)\xi^2/2}\,dx \,d\xi\Big|} \\  &\le &
\Big( \int_{-\infty}^\infty |\hat{f}(\xi)|^2 |\xi| ^{1/2}\, d \xi\Big)^{1/2} \Big(\int_{-\infty}^\infty \frac{d\xi}{|\xi|^{1/2}} \Big|\int_I e^{i(x\,\xi\, \sqrt{1+v(x)^2}- {v(x)}\xi^2/2)}\,dx\Big|^2\Big)^{1/2}.
\end{eqnarray*}
Here the first factor is controlled by the norm of $f$ in 
$W^{1/4}.$ Thus Lemma~\ref{essmaxest} will follow if we prove
that the second factor is bounded by some $C$. To this end,  we  write 
\begin{eqnarray}\label{osc}
\nonumber \lefteqn{\int_{-\infty}^\infty \frac{d\xi}{|\xi|^{1/2}} \Big|\int_I e^{i(x\,\xi\, \sqrt{1+v(x)^2}- {v(x)}\xi^2/2)}\,dx\Big|^2} \\&& = \int_{-\infty}^\infty \frac{d\xi}{|\xi|^{1/2}} \int \int_{I\times I} e^{i(x\,\xi\,  \sqrt{1+v(x)^2}- v(x)\xi^2/2)}e^{-i(y\,\xi\,\sqrt{1+v(y)^2}- {v(y)}\xi^2/2)}\, dx dy
 \nonumber \\ & &= \int \int_{I\times I}\, dx dy\, \int_{-\infty}^\infty \frac{e^{ia\,\xi  + ib\, \xi^2 }}{|\xi|^{1/2}}\,d\xi,
\end{eqnarray}
where $a = x\,\sqrt{1+v(x)^2}- y \sqrt{1+v(y)^2} $ and $b= (v(y)-v(x))/2.$ Observe that
\begin{multline*}
\big|\sqrt{1+v(x)^2} -  \sqrt{1+v(y)^2}\big| =  \frac{|v(x)+v(y)|\, |v(x)-v(y)|}{\sqrt{1+v(x)^2} +  \sqrt{1+v(y)^2}}  \le |v(x)-v(y)|.
\end{multline*}
In order to bound the last inner integral in (\ref{osc}), we shall distinguish between two cases.

\noindent Case 1:   $|y||b| <  |x-y|/4.$ Then we have 
\begin{eqnarray*} |a| &=& \left|(x-y) \sqrt{1+v(x)^2} + y \Big( \sqrt{1+v(x)^2}-\sqrt{1+v(y)^2}\Big)\right|\\
&>& |x-y| -\left|y \Big( \sqrt{1+v(x)^2}-\sqrt{1+v(y)^2}\Big)\right|\\ &\ge& |x-y| -|y ||v(x)-v(y)|\\ 
 &=& |x-y|-|y| |2b| \geq  |x-y|/2,
\end{eqnarray*}
and Lemma \ref{Peter} implies
\begin{eqnarray*}
\int_{-\infty}^\infty \frac{e^{ia\,\xi  + ib\, \xi^2 }}{|\xi|^{1/2}}\,d\xi
\le C|a|^{-\frac12} 
 \le C    
\frac1{|x-y|^{1/2}}.   
\end{eqnarray*}

\noindent Case 2:\quad $|y||b| \ge |x-y|/4$.    
   By using again Lemma \ref{Peter}, we conclude
$$\int_{-\infty}^{\infty} e^{ia\xi a + i b\xi^2 } \frac{d \xi}{|\xi|^{1/2}} \le C b^{-1/4} \le C \frac{|y|^{1/4}}{|x-y|^{1/4}}.$$

Summing up, we get for the iterated integral in (\ref{osc})
\[
\int \int_{I\times I}\, dx dy\, \int_{-\infty}^\infty \frac{e^{ia\,\xi  + ib\, \xi^2 }}{|\xi|^{1/2}}\,d\xi    \le  C
\int \int_{I\times I} \left(\frac1{|x-y|^{1/2}}  +
\frac{|y|^{1/4}}{|x-y|^{1/4}} \right)\,dxdy
\le C(I),
\]
and the proof of Lemma \ref{essmaxest} is complete.
\cqd

\vskip2mm

Next, we prove Theorem \ref{convergence}(ii). 
Because of Lemma \ref{fund}, it is sufficient to fix $s < 1/4$ and construct a 
$\varphi \in W^{s}$ for which the functions $ L_{iv/2} \varphi(x\,\sqrt{1+v^2}\,)$ diverge for $x$ in  a set of positive measure, as $v \rightarrow 0$
and $v$ avoids any given null set. 
The method is taken from \cite{DK}, though we prefer to make the construction
more explicit.

Choose a nonzero $f \in C^\infty_0$ supported in 
$\mathbb{R_-} = \{x: x<0\}$ and consider the functions $f_t(y) = f({y}/{t})e^{{2iy}/{t^2}}$ for small $t>0.$ Their Fourier transforms are $\widehat{f_t}(\xi) = t \hat{f}(t\, \xi-2/{t})$, and one finds that 
\begin{equation}\label{sob}
\|f_t\|_{W^s} \le C t^{1/2-2s} 
\end{equation}
for $t <1.$ Except for constant factors, $(1+v^2)^{1/4} L_{iv/2}f_t(x\,\sqrt{1+v^2}\,)$  is given by 
\begin{eqnarray}\label{1}
(1+v^2)^{1/4}v^{-1/2}\int_{\mathbb{R}} \exp\Big( i \frac1{2v}(x\, \sqrt{1+v^2}-y)^2\Big) f_t(y) dy.
\end{eqnarray}
Here we choose $\displaystyle v=v(x,t) = { x\, t^2}/{\sqrt{4-x^2\,t^4}}$
for $0<x<1$,  which implies 
\begin{equation}
  \label{eq:v}
   {v(x,t)}/{\sqrt{1+v(x,t)^2}} = {x\,t^2}/{2}.
\end{equation}

 Expanding the square in (\ref{1}) and using (\ref{eq:v}), we find that the expression (\ref{1}) for this value of $v$ and  $0<x<1$ equals $\sqrt{2} $ times 
\begin{multline*}
 \frac{1}{t\, \sqrt{x}}  \int_{\mathbb{R}}
\exp{\left(i\frac{2v(x,t)}{t^4}\right)}\exp{\left(-i\frac{2 y}{t^2}\right) }
\exp{\left(i \frac{y^2}{2v(x,t)}\right)} f(\,\frac{y}{t}\,)\, \exp{\left(2i\frac{y}{t^2}\right)}\, dy\\
=  \frac{1}{t\, \sqrt{x}} \exp{\left(i\frac{2v(x,t)}{t^4}\right)} \int_{\mathbb{R}}
\exp{\left(i \frac{y^2}{2v(x,t)}\right)} f(\,\frac{y}{t}\,)\, dy \\= 
 \exp{\left(i\frac{2v(x,t)}{t^4}\right)} \frac{1}{\sqrt{x}}\int_{\mathbb{R}} \exp{\Big(iy^2 \frac{\sqrt{1-{x^2t^4}/{4}}}{x}\Big)} f(y) dy \\ = \exp{\left(i\frac{2v(x,t)}{t^4}\right)} \frac{1}{\sqrt{x}} \Phi\Big( \frac{x}{\sqrt{1-{x^2t^4}/{4}}}\Big), 
 \end{multline*}
where 
$\displaystyle \Phi(z) = \int_{\mathbb{R}} f(y) \exp{\left(\frac{iy^2}{z}\right)} \,dy.$  This function  $\Phi$ is holomorphic in $\mathbb{C}\setminus\{0\}$
and  not identically $0$.  Thus there exists an interval $I \subset (1/2,1)$  such that $|\Phi(z)| > c$ for some constant $c >0$  when $z \in I.$ 
We can then find a subinterval $I' \subset I$ and an $\varepsilon > 0$ for which $x \in I'$ and $0 < t < \varepsilon$ imply  $\displaystyle {x}/{\sqrt{1- x^2 t^4/4}} \in I$ and thus $\displaystyle  |\Phi({x}/{\sqrt{1- x^2 t^4/4}}) | > c.$

To summarize the above, we have shown that  for some $c>0$
\begin{equation}\label{2}
\Big |\ L_{iv(x,t)/2} f_t\big(x \, \sqrt{1+v(x,t)^2}\,\big)\, \Big| > c,
\end{equation}
when $t< \varepsilon$ and $x \in I'.$
By continuity, one gets a stronger version of this inequality: it will
remain valid  if $v(x,t)$ is replaced by any
number in a sufficiently small neighborhood of  $v(x,t)$, a neighborhood
which may depend on $x$ and $t$.

We shall choose $\varphi = \sum_{j=1}^\infty j f_{t_j},$ where the numbers $t_j\in (0,\varepsilon)$ will be defined recursively. In particular, they shall satisfy $\sum_j j t_j^{1/2-2s} < \infty,$ which implies $\varphi \in W^s$
because of (\ref{sob}). Then
$$
 L_{iv/2} \varphi\big(x\, \sqrt{1+v^2}\,\big) = \sum_{j=1}^\infty j L_{iv/2} f_{t_j}\big(x\,\sqrt{1+v^2}\,\big).
$$
Now consider  $x\in I'$ and any $k=1,2, \dots$.
Our idea is  to make sure that for   $v$ close to
$v(x,t_k)$,
the term with $j=k$  is dominating in the above sum. More precisely, we shall have 
\begin{eqnarray}\label{3}
\Big |\, L_{iv/2} f_{t_j} ( x \, \sqrt{1+v^2}\,) \, \Big| < 2^{-j}, \qquad j\neq  k, 
\end{eqnarray}
 for   $x \in I'$ and $1/2 < v/v(x,t_k) < 2$. Combining this with (\ref{2})
and its stronger version,  we see that  for   $x \in I'$ and $v $ close to
$v(x,t_k)$, 
$$ \Big |\, L_{iv/2}\varphi\big(x\, \sqrt{1+v^2 }\, \big)\, \Big| \ge ck -\sum_{j\neq k} j\, 2^{-j}.$$
The right-hand side here
 tends to $+\infty$ with $k$, and divergence will follow    once we have established (\ref{3}).

In the recursive construction of the $t_j,$ we  start with any $t_1 \in (0,\varepsilon).$ Assume now $t_1,\dots, t_{J-1}$ chosen so that (\ref{3}) holds when $j,k < J.$ Then we must find $t_J$ so that, when $x \in I'$,
\begin{equation}\label{4}
\Big |\ L_{iv/2}f_{t_j}\big(x\, \sqrt{1+v^2}\,\big)\, \Big| < 2^{-j}. \,\qquad j=1, \dots ,J-1,
\end{equation}
for $1/2 < v/v(x,t_J) < 2$, and
\begin{equation}\label{5}
\Big | L_{iv/2}f_{t_J}\big(x\, \sqrt{1+v^2}\,\big) \Big| < 2^{-J} \quad 
\mbox{for} \quad  
  \frac12 < v/v(x,t_k) < 2, \quad   k=1, \dots ,J-1.
\end{equation}
Aiming at (\ref{4}), we observe that each $f_{t_j}$ is a $C^\infty_0$ function and so 
$L_{is} f_{t_j} \rightarrow f_{t_j}$ uniformly in $\mathbb{R}$ as $s \rightarrow 0^+.$  Now $v(x,t) \rightarrow 0$ as $t \rightarrow 0,$ and $I' \subset (1/2,1)$ but the $f_{t_j}$ are supported in $\mathbb{R}_-$. This means that (\ref{4}) will hold for the indicated values of $x$ and $v$, if $t_J$ is chosen small enough.

To obtain (\ref{5}), we simply estimate $L_{iv/2}f_{t_J}$ by the supremum norm of the kernel  $L_{{iv/2}}$ times the $L^1$ norm of $f_{t_J}.$ This  product is 
$Cv^{-1/2}t_J,$ and (\ref{5}) follows if $t_J$ is small.
The recursive construction and the proof of Theorem~\ref{convergence}(ii)
are complete.

\vskip2mm

\noindent \textit{Proof of Theorem \ref{convergence}(iii).} \hskip2mm
The analog of Part (i) for $W^s_H$ is obvious, since  
$W^s_H \subset W^s$; see \cite[Theorem 3(i)]{BT}. As for Part (ii),
 observe that  the function $\varphi$ constructed
above is in $W^s$ and has compact support. But then $\varphi$ is also  
in $W_H^s$, as proved in \cite[Theorem 3(iii)]{BT}.

Theorem \ref{convergence} is completely proved.
\cqd

\section{Proof of theorem \ref{noway}}

Lemma \ref{fund}
implies that $\mathcal{M}f(x)$ can be estimated from below by a positive 
constant times
\begin{equation}
  \label{sup}
\underset{0<v < 1}{\mbox{\textrm{ess sup}}}\;|L_{iv/2}f(x\,\sqrt{1+v^2}\,)|.
\end{equation} 

We first consider the case $p<\infty$.
Fix a large $x_0 >0$ and choose a function  $0 \le \tau \in C_0^\infty,$ with  ${\rm supp  }\,\, \tau \subset (-1,1).$ Let  $f$ be given  by $\hat{f}(\xi) = 2\pi e^{-i\,x_0 \, \xi}\tau(\xi),$ and  define $v(x) \in (0,1)$ by $x\,\sqrt{1+v(x)^2} = x_0$ for $x_0/\sqrt{2} < x  < x_0$. Then for these $x,$  
\begin{multline*}
L_{iv(x)/2}f(x\,\sqrt{1+v(x)^2}\,) = \frac1{2\pi} \int_{\mathbb{R}} e^{-iv(x)\xi^2/2} e^{i\, x\, \xi\sqrt{1+v(x)^2}} \hat{f}(\xi)\, d \xi \\ = \int_{\mathbb{R}} e^{-i v(x) \xi^2/2}\, e^{i x_0 \xi} \, e^{-i x_0 \xi}\, \tau(\xi)\, d\xi = \int e^{-iv(x) \xi^2/2} \tau(\xi) \,d \xi.
\end{multline*}
For $\xi \in \; \mbox{supp}\, \tau$ one  has $0 < v(x) \xi^2/2 < 1/2,$ and so 
\begin{equation}
  \label{eq:P}
\Re ( L_{i v(x) /2} f(x\, \sqrt{1+v(x)^2}\,)) > \cos \frac{1}{2}\, \int \tau >0,
\qquad x_0/\sqrt{2} < x < x_0.
  \end{equation}
 By continuity, this holds also if
the value of $v(x)$ is slightly modified. Thus  
$\| \mathcal{M}f\|_p \ge cx_0^{1/{p}}$ for some $c>0,$ 
and  the weak $L^p$ quasinorm of $\mathcal{M}f$ satisfies the same
inequality.  But 
$$ \|f\|^2_{W^s} = \int |\hat{f}(\xi)|^2(1+|\xi|^2)^{s} d\xi = 4\pi^2 \int \tau(\xi)^2(1+|\xi|^2)^{s} \,d \xi$$
is independent of $x_0.$ Finally let $x_0 \rightarrow +\infty,$ 
to get the desired unboundedness. 

For $p=\infty$ we first assume that $s>1/2$. H\"older's inequality
then implies that $\|\hat{f}\|_{L^1} \le C \|f\|_{W^s}$.
Thus for any $x$ and any $v$ one can estimate
\[
L_{iv/2}f(x\,\sqrt{1+v^2}\,) = \frac1{2\pi} \int_{\mathbb{R}} e^{-iv\xi^2/2} e^{i\, x\, \xi\sqrt{1+v^2}} \hat{f}(\xi)\, d \xi
\]
by means of the $W^s$ norm of $f$, as required.

To find a counterexample for $p=\infty$ and  $s\le 1/2$, we
 modify the above construction by taking now
 $0 \le \tau \in C^\infty,$ supported in $\RR_+$ and such that
$\tau(\xi) = \xi^{-1} (\log \xi)^{-2/3}$ for $\xi > 2$. As before,
 $\hat{f}(\xi) = 2\pi e^{-i\,x_0 \, \xi}\tau(\xi),$ 
but  $x_0 > 0$ is now fixed. One easily verifies that $f\in W^s$.
The choice of $v(x)$ is again given by
$x\,\sqrt{1+v(x)^2} = x_0$, but now only when 
$x$ is in the interval 
\[
I = \left(\frac{x_0}{\sqrt{1+v_0^2}},\; \frac{x_0}{\sqrt{1+v_0^2/4}}\right),
\]
 for
some small $v_0$.  Then  $v_0/2<v(x)<v_0$, and for almost all  $x \in I$ 
we conclude essentially as before that
\begin{equation*}
L_{iv(x)/2}f(x\,\sqrt{1+v(x)^2}\,)  = \int e^{-v(x) \xi^2/2} \tau(\xi)\, d \xi.
\end{equation*}
Notice that since this is now obtained via a truncation of 
$f$ at $+\infty$ and
an $L^2$ limit, the integral here should be evaluated as  
$\lim_{R\to +\infty} \int_0^R$. 
Part of this integral can be estimated as above; indeed
\[
\Re\left( \int_0^{1/\sqrt{v_0}} e^{-iv(x) \xi^2/2} \tau(\xi) \,d\xi\right)
>  \cos \frac{1}{2}\, \int_0^{1/\sqrt{v_0}} \tau(\xi) \,d \xi 
\ge c \left(\log \frac1{\sqrt{v_0}}\right)^{1/3}.
\]
In the remaining part, we integrate by parts and get
\[
\left| \int_{1/\sqrt{v_0}}^\infty e^{-iv(x) \xi^2/2} \tau(\xi) \,d\xi\right|
\le \frac1{v(x)} \frac{\tau(1/\sqrt{v_0})}{1/\sqrt{v_0}}
+ \frac1{v(x)} \int_{1/\sqrt{v_0}}^\infty \left|\frac{d}{d\xi}
\frac{\tau(\xi)}{\xi}\right|\,d\xi.
 \]
The last integral  equals ${\tau(1/\sqrt{v_0})}{\sqrt{v_0}}$,
because the derivative in the integrand is  negative here.
Since $v(x) > v_0/2$, each term in the above right-hand side is at most
 $2\log(1/\sqrt{v_0})^{-2/3}$.

Summing up, we see that
\[
\big|L_{iv(x)/2}f(x\,\sqrt{1+v(x)^2}\,)\big| \ge c \left(\log \frac1{v_0}\right)^{1/3}
\]
for a.a. $x \in I$, also after a slight
modification of $v(x)$.
Letting $v_0 \to 0$, we conclude that the essential supremum in (\ref{sup})
is not in $L^\infty$ for this $f$,
which ends the proof.

 \cqd

\end{document}